\documentclass[twocolumn]{article}
\usepackage{epsfig}
\usepackage{graphicx}
\usepackage{color}
\usepackage[acmtitlespace,nocopyright]{acmneww}
\usepackage{times}
\usepackage{amsmath}
\usepackage{amssymb}
\usepackage{xspace}
\usepackage{txfonts}
\usepackage{enumitem}

\newcommand {\mm}[1] {\ifmmode{#1}\else{\mbox{\(#1\)}}\fi}

\newcommand {\scalprod}[2] {{\langle #1 , #2 \rangle}}
\newcommand {\denselist}{\itemsep 0pt\parsep=1pt\partopsep 0pt}

\newcommand{\proof}{\noindent{\sc Proof.~}}

\newcommand{\eop}{\hfill\usebox{\smallProofsym}\bigskip}  %
\newsavebox{\smallProofsym}                            
\savebox{\smallProofsym}                               %
{
\begin{picture}(6,6)                                   %
\put(0,0){\framebox(6,6){}}                            %
\put(0,2){\framebox(4,4){}}                            %
\end{picture}                                          %
}                                                      %
\makeatletter
\long\def\@makecaption#1#2{%
  \vskip\abovecaptionskip
  \sbox\@tempboxa{\small #1: #2}%
  \ifdim \wd\@tempboxa >\hsize
    \small #1: #2\par
  \else
    \global \@minipagefalse
    \hb@xt@\hsize{\hfil\box\@tempboxa\hfil}%
  \fi
  \vskip\belowcaptionskip}
\makeatother

\newcommand{\Rspace}        {\mm{{\mathbb R}}}
\newcommand{\Dtri}          {\mm{{\mathcal D}}}
\newcommand{\Ktri}          {\mm{{\mathcal K}}}

\newcommand{\Delaunayf}[1]  {\mm{\rm Df}{({#1})}}
\newcommand{\Rajanf}[1]     {\mm{\rm Rf}{({#1})}}
\newcommand{\Voronoif}[1]   {\mm{\rm Vf}{({#1})}}
\newcommand{\Radiusf}[1]    {\mm{\rm Rad}_{#1}}
\newcommand{\N}[2]          {\mm{N_{#1}{({#2})}}}
\newcommand{\NV}[2]         {\mm{NV_{#1}{({#2})}}}
\newcommand{\Sd}[1]         {\mm{\rm Sd\,}{#1}}
\newcommand{\barycenter}[1] {\mm{\rm b}{({#1})}}
\newcommand{\circumcenter}[1] {\mm{\rm c}{({#1})}}
\newcommand{\circumradius}[1] {\mm{\rm r}{({#1})}}
\newcommand{\orientation}[1]  {\mm{\rm sgn}{({#1})}}
\newcommand{\Height}        {\mm{\rm H}}
\newcommand{\Star}[1]       {\mm{\rm St\,}{#1}}
\newcommand{\us}[1]         {\mm{|{#1}|}}
\newcommand{\ksx}           {\mm{\kappa}}

\newcommand{\volume}[1]     {\mm{\rm vol}{({#1})}}
\newcommand{\area}[1]       {\mm{\rm area}{({#1})}}
\newcommand{\conv}[1]       {\mm{\rm conv\,}{#1}}
\newcommand{\dist}[2]       {\mm{\|{#1}-{#2}\|}}
\newcommand{\norm}[1]       {\mm{\|{#1}\|}}
\newcommand{\diff}          {\mm{\rm \,d}}

\newtheorem{result}{}

\newcommand{\ignore}[1]{}

\title{The Voronoi Functional is Maximized ~\\
       by the Delaunay Triangulation in the Plane
       \thanks{This research is partially supported by
               the Russian Government under the Mega Project 11.G34.31.0053,
               by the {\sc Toposys} project FP7-ICT-318493-STREP,
               by ESF under the ACAT Research Network Programme,
               by RFBR grant 11-01-00735, and by NSF grants DMS-1101688, DMS-1400876.}
       }

\author{Herbert Edelsbrunner\thanks{IST Austria (Institute of Science and
          Technology Austria), Kloster\-neu\-burg, Austria.},
        Alexey Glazyrin\thanks{Department of Mathematics,
         The University of Texas at Brownsville, Texas, USA.},
        Oleg R.\ Musin${}^{\ddagger}$ and
        Anton Nikitenko${}^{\dagger}$
}

\begin{document}
\maketitle

\begin{abstract}
  We introduce the Voronoi functional of a triangulation of a finite
  set of points in the Euclidean plane and prove that
  among all geometric triangulations of the point set,
  the Delaunay triangulation maximizes the functional.
  This result neither extends to topological triangulations in the plane
  nor to geometric triangulations in three and higher dimensions.
\end{abstract}

\vspace{0.1in}

\section{Introduction}
\label{sec1}

The \emph{Voronoi diagram} of points $x_1, x_2, \ldots, x_m$
in $n$-dimen\-sion\-al Euclidean space decomposes $\Rspace^n$
into $m$ convex polyhedra, called \emph{Voronoi domains}.
The domain of $x_i$ is the set of points $x \in \Rspace^n$
for which $x_i$ minimizes the Euclidean distance:
\begin{align}
  V_i  &=  \{ x \in \Rspace^n
    \mid \dist{x}{x_i} \leq \dist{x}{x_j}, 1 \leq j \leq m \}.
\end{align}
For points in general position, the geometrically realized nerve
(the straight-line dual in $\Rspace^2$)
decomposes the convex hull of the points into simplices,
called the \emph{Delaunay triangulation}.
It is a particular \emph{geometric triangulation} of the points,
which is a simplicial complex whose vertices are the given points and whose
underlying space is convex.
Delaunay triangulations are used in numerous applications
and often preferred over other geometric triangulations.
The question arises why the Delaunay triangulation is better than others.
In the plane, the advantages of Delaunay triangulation are
sometimes rationalized by the max-min angle criterion \cite{Sib78}.
It requires that the diagonal of every convex quadrangle
should be chosen to maximize the minimum
of the six angles in the two triangles making up the quadrangle
\cite{Law77}.
It follows that the Delaunay triangulation lexicographically maximizes
the non-decreasing sequence of angles.

Triangulations of the same finite collection of points in
$n \geq 3$ dimensions can have different length sequences of angles.
Instead of comparing them directly, we may consider a functional,
which assigns a real number to every triangulation.
In this paper, we consider functionals that are defined for $n$-simplices
and we assign the sum of the values over all $n$-simplices to the triangulation.
For example, instead of the angle sequence we may consider
the functional that maps a triangulation to the sum of minimum angles
of $n$-simplices.
In $\Rspace^2$, this functional attains
its minimum for the Delaunay triangulation.
Other functionals that attain their minimum for the Delaunay triangulation
in the plane have been studied in \cite{CXGL10,Lam94,Mus97,Mus10,Raj94,Rip90}
and generalized from finite collections to (infinite) Delone sets
in \cite{DEGM13}.
For example, the Rajan functional is defined by mapping
a triangle $\Delta$ with edges of length $a, b, c$
to $a^2 + b^2 + c^2$ times the area of $\Delta$.
There is an intuitive geometric interpretation obtained by
lifting a point with coordinates $\xi$ and $\eta$
to $(\xi, \eta, \xi^2 + \eta^2)$,
which is a point on the standard paraboloid in $\Rspace^3$.
Similarly, we lift $\Delta$ to a triangle in space,
namely the one spanned by the three lifted vertices.
The Rajan functional maps $\Delta$ to the volume
between the paraboloid and the lifted triangle,
of course restricted to the vertical prism over $\Delta$.
The Rajan functional of a geometric triangulation is therefore the
volume between the paraboloid and the lifted triangulation.
Since the Delaunay triangulation corresponds to taking the convex hull
of the lifted points \cite{Ede01},
it is easy to see that the Delaunay triangulation minimizes this functional.

The focus in this paper is a related concept, which we refer to
as the \emph{Voronoi functional}; see also \cite{Ako09}.
For a triangle $\Delta$ with acute angles, it is the volume
below the standard paraboloid and above the tangent planes that
touch the paraboloid in the lifted vertices of $\Delta$,
again restricted to the vertical prism above $\Delta$.
The definition can be extended to triangles with non-acute angles
and to $n$-simplices, as we will explain in the body of this paper.
Our main results are a characterization of the Voronoi functional
for triangulations in $\Rspace^2$ in terms of squared Euclidean
distances to the given points,
and a proof that the Voronoi functional attains its maximum
for the Delaunay triangulation;
see the Voronoi Cell Decomposition
and the Voronoi Optimality Theorems in Section \ref{sec3}.
While the Voronoi functional is defined for topological triangulations
of the given points, which allow for folding along edges,
the optimality of the Delaunay triangulation holds only within
the smaller class of geometric triangulations.
Similarly, we have a counterexample to the optimality of the
Delaunay triangulation among the geometric triangulations
of points in three and higher dimensions.

\paragraph{Outline.}
Section \ref{sec2} defines the Voronoi functional
for geometric and topological triangulations of points in $\Rspace^2$.
Section \ref{sec3} introduces the circumcenter map, which is instrumental
for the interpretation of the Voronoi functional as a volume in $\Rspace^3$,
and it proves our two main results.
Section \ref{sec4} gives counterexamples to extending the optimality result.
Section \ref{sec5} concludes the paper.

\section{The Voronoi Functional}
\label{sec2}

We introduce the Voronoi functional in three steps,
beginning with the relatively easy case of an acute triangle.

\begin{figure}[hbt]
 \centering \resizebox{!}{2.4in}{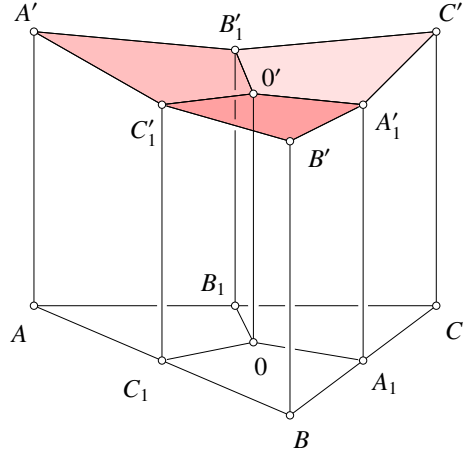}
 \caption{Lifting the vertices of an acute triangle,
   we are interested in the volume below the paraboloid and above the
   three planes.}
 \label{fig:Voronoi-acute}
\end{figure}

\paragraph{Acute case.}
Let $A, B, C$ be the vertices of an acute triangle, $\Delta$, in $\Rspace^2$.
For a point $x$, we write $\N{\Delta}{x}$ for the vertex among the three
that is nearest to $x$, and we define
\begin{align}
  \Voronoif{\Delta}  &=  \int_{x \in \Delta} \dist{x}{\N{\Delta}{x}}^2 \diff x .
  \label{eqn:Voronoi-1}
\end{align}
To interpret $\Voronoif{\Delta}$ geometrically, we introduce the
unit paraboloid as the graph of
$\varpi \colon \Rspace^2 \to \Rspace$ defined by $\varpi (x) = \norm{x}^2$.
Setting $f_A (x) = 2 \scalprod{x}{A} - \norm{A}^2$,
we note that the graph of $f_A \colon \Rspace^2 \to \Rspace$
is the tangent plane that touches the paraboloid in the point
$A' = (A, \norm{A}^2)$.
The vertical distance between the paraboloid and the plane above
a point $x \in \Rspace^2$ is
\begin{align}
  \varpi(x) - f_A(x)  &=  \norm{x}^2 - 2 \scalprod{x}{A} + \norm{A}^2 
                       =  \dist{x}{A}^2 .
\end{align}
Subdividing $\Delta$ into regions of constant nearest vertex,
we get three quadrangles.
Integrating the squared distance to the vertex over each quadrangle,
and adding the results, we get $\Voronoif{\Delta}$ as the volume
between the paraboloid and the upper envelope of the three tangent planes;
see Figure \ref{fig:Voronoi-acute}.

Before extending $\Voronoif{\Delta}$ to obtuse triangles,
we develop an algebraic formula in the acute case.
We use the notation of Figure \ref{fig:Voronoi-acute},
assuming $\Rspace^2$ is spanned by the first two coordinate axes
of $\Rspace^3$, and $0 = (0,0,0)$
is the common origin of $\Rspace^2$ and $\Rspace^3$.
Without loss of generality, we may assume that $0$ is the
circumcenter of $\Delta$.
Let $\area{\Delta}$ be the area of the triangle, and
$\volume{\Delta}$ the volume of the triangular prism
bounded from below by $\Delta$ and from above by the triangle
with vertices $A', B', C'$.
Letting $R$ be the radius of the circumcircle, we have
$\norm{A}^2 = \norm{B}^2 = \norm{C}^2 = R^2$, and therefore
$\volume{\Delta}  =  R^2 \cdot \area{\Delta}$.
Writing $Q_A = A C_1 0 B_1$, $Q_B = B A_1 0 C_1$, and $Q_C = C B_1 0 A_1$
for the quadrangles subdividing the acute triangle,
we let $\volume{Q_A}$ be the volume of the quadrangular prism between
$Q_A$ and $A' C_1' 0' B_1'$,
and similar for $Q_B$ and $Q_C$.
Unlike suggested by Figure \ref{fig:Voronoi-acute},
the point $0'$ is below $0$, so that a portion of each quadrangular prism
is below the horizontal coordinate plane.
The notion of volume we use is signed,
which means that $\volume{Q_A}$ is the (unsigned) volume of the portion
above $\Rspace^2$ minus the (unsigned) volume of the portion below $\Rspace^2$.
Writing $a = \dist{B}{C}$, $b = \dist{C}{A}$, and $c = \dist{A}{B}$
for the lengths of the three edges, we set
\begin{align}
  \Rajanf{\Delta}  &=  \tfrac{\area{\Delta}}{12} \cdot (a^2 + b^2 + c^2) .
\label{rajan}
\end{align}
This is known as the \emph{Rajan functional},
which is the volume between the paraboloid and the triangle $A' B' C'$.
It is now easy to express the Voronoi functional by subtracting the
Rajan functional and volumes of the three quadrangular prisms
from the volume of the triangular prism:
$\Voronoif{\Delta} = \volume{\Delta} - \Rajanf{\Delta}
                  - \volume{Q_A} - \volume{Q_B} - \volume{Q_C}$.
After half a page of trigonometric calculations, we get
\begin{align}
  \Voronoif{\Delta}  &=  \tfrac{\area{\Delta}}{12} \cdot (a^2 + b^2 + c^2 - 4 R^2) ;
  \label{eqn:Voronoi-2}
\end{align}
see Appendix B.

\begin{figure}[hbt]
 \centering \resizebox{!}{1.5in}{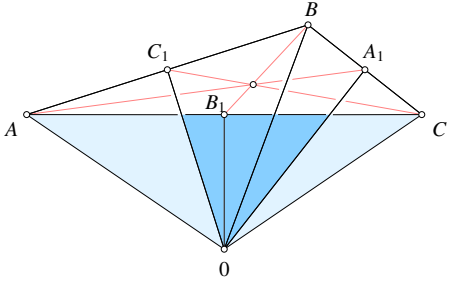}
 \caption{Two of the six triangles in the subdivision of the
   triangle spanned by $A, B, C$ have negative orientation.
   In the blue shaded region, points contribute two terms to
   the Voronoi functional,
   and in the dark blue region, these two terms do not cancel.}
 \label{fig:fold}
\end{figure}

\paragraph{Obtuse case.}
To extend the formula to the case of an obtuse triangle, it is convenient
to further subdivide each quadrangle into two triangles by drawing the
edges between $0$ and the three vertices.
Recalling that $B_1$ is the midpoint of the edge from $A$ to $C$, we set
\begin{align}
  \mu(A B_1 0)  &=  \int_{x \in A B_1 0} \dist{x}{A}^2 \diff x ,
\end{align}
and similarly for the other five triangles in the subdivision.
In the acute case, we can rewrite \eqref{eqn:Voronoi-1} to get
\begin{align}
  \Voronoif{\Delta}  &= \mu(A 0 B_1) + \mu(A C_1 0) + \mu(B 0 C_1) \\
                    &+ \mu(B A_1 0) + \mu(C 0 A_1) + \mu(C B_1 0) .
\end{align}
Let now $\Delta$ be obtuse and use \eqref{eqn:Voronoi-2}
to define the Voronoi functional of $\Delta$.
Assuming the angle at $B$ exceeds $90^\circ$, we get
\begin{align}
  \Voronoif{\Delta}  = &- \mu(A 0 B_1) + \mu(A C_1 0) + \mu(B 0 C_1) \\
                      &+ \mu(B A_1 0) + \mu(C 0 A_1) - \mu(C B_1 0) ;
\end{align}
see Appendix A.
This has a geometric interpretation,
which we illustrate in Figure \ref{fig:fold}.
In particular, two of the six triangles subdividing $\Delta$
have negative orientation, namely $A B_1 0$ and $C 0 B_1$,
and we record the volume of the corresponding triangular prisms
with a minus sign.
The reason for the negative orientation is that $0$ lies outside the triangle
and, in the illustrated case, on the opposite side of the edge from $A$ to $C$.

Besides the fold-over interpretation of $\Voronoif{\Delta}$,
we will use the interpretation as the integral over a difference between
two squared distances.
Call a vertex \emph{visible} from $x$ if the line segment connecting
$x$ to the vertex does not intersect $\Delta$ other than in this vertex.
Writing $\NV{\Delta}{x}$ for the nearest vertex of $\Delta$ that is
visible from $x$, we have
\begin{align}
  \Voronoif{\Delta}  &=  \int_{x \in \Rspace^2}
    \left( \dist{x}{\N{\Delta}{x}}^2
         - \dist{x}{\NV{\Delta}{x}}^2 \right) \diff x .
  \label{eqn:Voronoi-3}
\end{align}
For points $x$ inside $\Delta$,
$\NV{\Delta}{x}$ is not defined and we set $\dist{x}{\NV{\Delta}{x}}^2 = 0$.
This way we get a formula that is correct both for acute and
for obtuse triangles.
Indeed, in the acute case, we have $\N{\Delta}{x} = \NV{\Delta}{x}$
for all points $x$ outside $\Delta$,
which implies that \eqref{eqn:Voronoi-3} agrees with \eqref{eqn:Voronoi-1}.
In the obtuse case, we have $\N{\Delta}{x} \neq \NV{\Delta}{x}$ for all points
$x$ in the shaded region shown in Figure \ref{fig:fold}.
In this case, $\Voronoif{\Delta}$ is not the volume between the paraboloid
and the three tangent planes restricted to the triangular prism,
but rather something smaller than this volume.

\paragraph{In the large.}
We extend the Voronoi functional by taking the
sum over all triangles of a triangulation.
Specifically, letting $\Ktri$ be a triangulation of $S \subseteq \Rspace^2$,
we set
\begin{align}
  \Voronoif{\Ktri}  &=  \sum_{\Delta \in \Ktri} \Voronoif{\Delta} .
  \label{eqn:Voronoi-4}
\end{align}
What do we mean by a triangulation of $S$?
Most important is the
\emph{Delaunay triangulation} that consists of all triangles spanned
by three points of $S$ such that no point of $S$ is enclosed
by the circumcircle of the triangle.
Assuming the points are in general position -- by which we mean
that no four points lie on a common circle --
the Delaunay triangulation is well defined and unique.

A more general notion is a \emph{geometric triangulation},
which is a simplicial complex in $\Rspace^2$ whose vertex set is $S$
and whose underlying space is $\conv{S}$.
Clearly, the Delaunay triangulation is a geometric triangulation of $S$,
and in many ways, it is the most natural and most interesting
geometric triangulation of the point set.

More general yet is a \emph{topological triangulation},
which is a simplicial complex homeomorphic to a disk whose vertex set is $S$.
In contrast to a geometric triangulation,
the triangles of a topological triangulation may intersect as we think
of them abstractly, worrying primarily about how they are connected.
Every geometric triangulation is also a topological triangulation
but not the other way round.
If no three points of $S$ are collinear,
then every triangle of a topological triangulation maps to a geometric
triangle, whose orientation may be positive or negative.
We interpret \eqref{eqn:Voronoi-4} accordingly,
namely that the sign of $\Voronoif{\Delta}$ is the same as that of
the orientation of the triangle.

\section{Optimality}
\label{sec3}

We give a complete description of the Voronoi functional for Delaunay triangulations,
and we prove that among the geometric triangulations of a finite set in $\Rspace^2$,
the Delaunay triangulation maximizes the Voronoi functional.

\begin{figure}[hbt]
 \centering \resizebox{!}{1.8in}{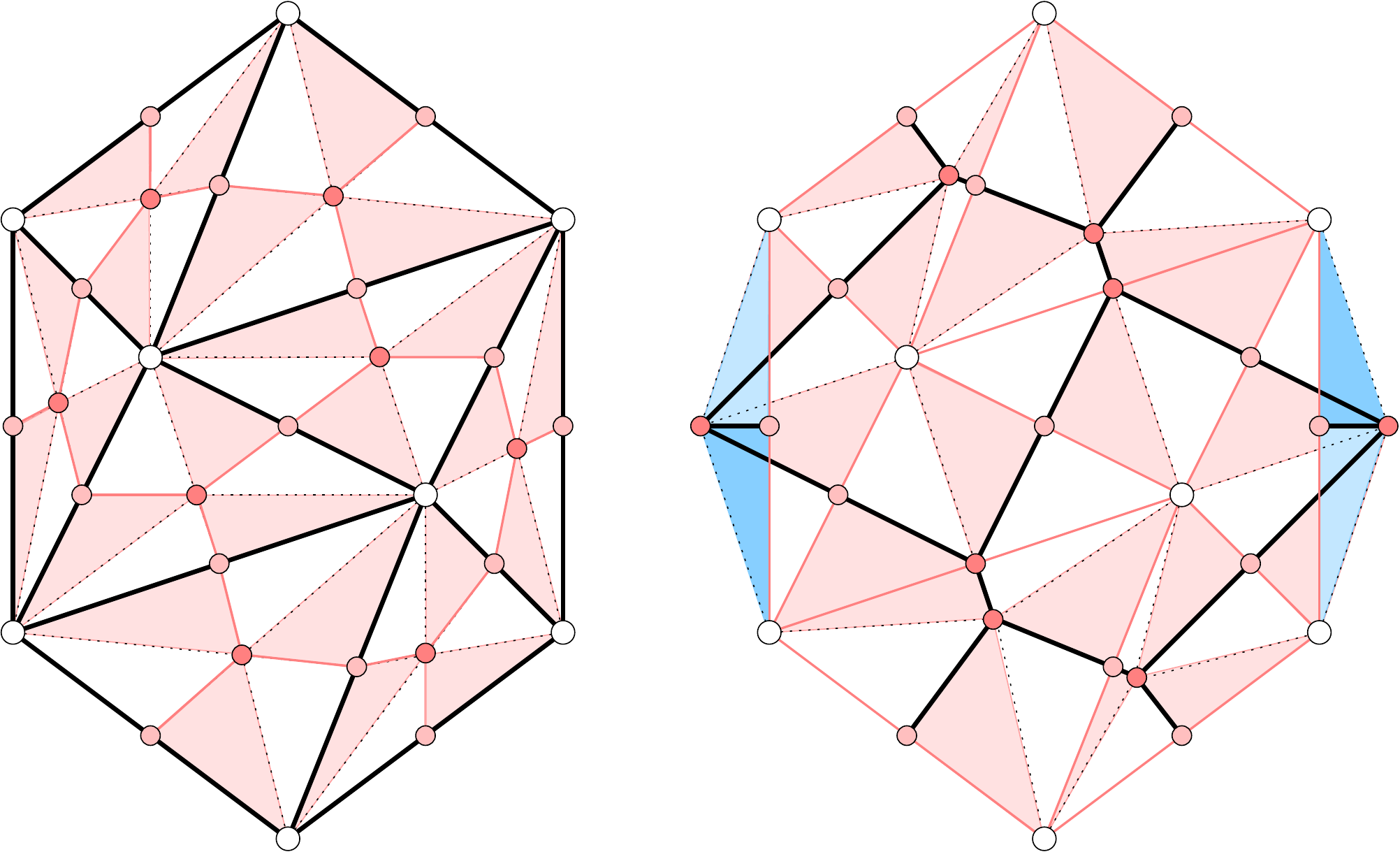}
 \caption{\emph{Left:} the barycentric subdivision of a Delaunay
   triangulation with triangles shaded in checkerboard style.
   \emph{Right:} the image of the barycentric subdivision under
   the circumcenter map.
   Four triangles flip over, and four triangles are squeezed to line segments.}
 \label{fig:Sd}
\end{figure}

\paragraph{The circumcenter map.}
Let $\Ktri$ be a geometric triangulation of $S \subseteq \Rspace^2$.
The \emph{barycenter} of a simplex $\ksx \in \Ktri$ is the point
$\barycenter{\ksx} \in \Rspace^2$
that is the average of the vertices of $\ksx$.
A \emph{flag} a sequence of simplices in $\Ktri$ such that each
simplex is a proper face of its successor.
The \emph{barycentric subdivision} of $\Ktri$ is the simplicial complex,
$\Sd{\Ktri}$, whose vertices are the barycenters of the simplices in $\Ktri$
and whose simplices correspond to flags in $\Ktri$; see Figure \ref{fig:Sd}.
Every $\ksx \in \Ktri$ has a unique smallest circle that passes through
its vertices.
We call the center of this circle the \emph{circumcenter},
$\circumcenter{\ksx}$, and its radius the \emph{circumradius},
$\circumradius{\ksx}$, of the simplex.
We introduce the \emph{circumcenter map} and the \emph{height map},
\begin{align}
  \Gamma  &\colon \us{\Sd{\Ktri}} \to \Rspace^2 , \\
  \Height &\colon \us{\Sd{\Ktri}} \to \Rspace ,
\end{align}
each the piecewise linear extension of a vertex map:
the first such that $\Gamma (\barycenter{\ksx}) = \circumcenter{\ksx}$,
and the second such that $\Height (\barycenter{\ksx})
= \norm{\circumcenter{\ksx}}^2 - \circumradius{\ksx}^2$,
for every $\ksx \in \Ktri$.
As illustrated in Figure \ref{fig:Sd}, the circumcenter map distorts the
barycentric subdivision so it aligns with the Voronoi diagram,
with the exception of occasional fold-edges caused by obtuse triangles.
The height map lifts the image of the circumcenter map to $\Rspace^3$.
To see how it does it, let $A \in S$ be a vertex of $\Ktri$,
and consider its star in $\Sd{\Ktri}$.
A triangle, $abc$, in this star is spanned by the barycenters
$a = \barycenter{A}$, $b = \barycenter{AB}$, $c = \barycenter{ABC}$,
with $AB$ an edge, and $ABC$ a triangle in $\Ktri$.
The images under $\Gamma$ are the point, $\Gamma (a) = a$,
the midpoint of the edge, $\Gamma (b) = b$,
and the circumcenter of the triangle, $\Gamma (c) = \circumcenter{ABC}$.
The corresponding heights are
$\Height (a) = f_A ( \Gamma (a))$,
$\Height (b) = f_A ( \Gamma (b))$, and
$\Height (c) = f_A ( \Gamma (c))$.
We see that the combined map, $( \Gamma , \Height )$,
sends the entire star of $A$ to the plane that is the graph of $f_A$.
Similarly, the star of $B$ is sent to the graph of $f_B$,
and the star of $C$ is sent to the graph of $f_C$.
The stars are glued along shared boundary pieces,
which must therefore lie in the common intersection of the planes.
Indeed, we have $f_A (\Gamma (b)) = f_B (\Gamma (b))$ and
$f_A (\Gamma (c)) = f_B (\Gamma (c)) = f_C (\Gamma (c))$.

We use the two maps to recast the Voronoi functional as an integral.
Write $\det (\delta)$ and $\det ( \Gamma (\delta))$ for the signed
areas of a triangle and its image under the circumcenter map.
With this, we have
\begin{align}
  \Voronoif{\delta}  &=  \frac{\det(\Gamma(\delta))}{\det(\delta)} \cdot
    \int_{x \in \delta} \left( \norm{\Gamma (x)}^2 - \Height (x) \right) \diff x .
  \label{eqn:detoverdet}
\end{align}
The ratio in front of the integral captures the area distortion
experienced by $\delta$, and it is negative iff $\Gamma$ reverses the
orientation of the triangle.
Geometrically, we interpret $\Voronoif{\delta}$ as the
signed volume between the paraboloid and the lifted copy of $\Gamma (\delta)$.
Finally, $\Voronoif{\Ktri}$ is the sum of the $\Voronoif{\delta}$,
over all triangles $\delta \in \Sd{\Ktri}$.

\paragraph{Inclusion-exclusion.}
In the case in which $\Ktri = \Dtri$ is the Delaunay triangulation
of $S \subseteq \Rspace^2$ and all angles are acute,
$\Voronoif{\Dtri}$ has a very appealing geometric interpretation
as the volume of the body between the paraboloid and the upper
envelope of the graphs of $f_A$, $A \in S$,
restricted to within the vertical prism over the convex hull of $S$.
This interpretation applies more generally,
namely as long as the angles opposite to the convex hull edges are acute.

\begin{figure}[hbt]
 \centering \resizebox{!}{1.6in}{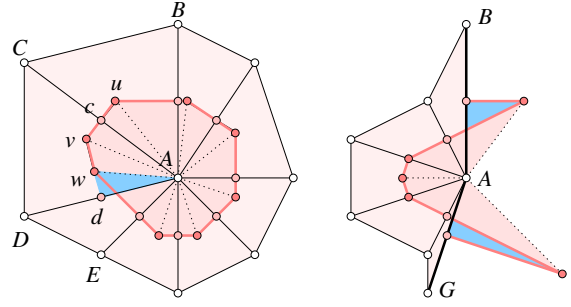}
 \caption{Images of stars in the barycentric subdivision of
   the Delaunay triangulation.
   An interior vertex on the \emph{left},
   and a boundary vertex with obtuse angles opposite the boundary edges
   on the \emph{right}.}
 \label{fig:stars}
\end{figure}
To explain why this is true, we consider an interior vertex $A$;
see Figure \ref{fig:stars} on the left.
Let $ABC$ and $ACD$ be two triangles and $uv$ the corresponding edge
of the Voronoi polygon.
Letting $c$ be the midpoint of the shared edge, $AC$,
we have $uc$ and $cv$ as images under $\Gamma$ of two edges in $\Sd{\Ktri}$.
If the angles at $B$ and $D$ inside the mentioned triangles are acute,
then $c$ lies in the interior of the Voronoi edge,
the triangles $Auc$ and $Acv$ both have the same orientation as
their preimages under $\Gamma$, and they decompose the triangle $Auv$.
To describe the other case, we consider the triangles $ACD$ and $ADE$
with corresponding edge $vw$ of the Voronoi polygon in Figure \ref{fig:stars}.
If the angle at $E$ is obtuse, then the angle at $C$ is necessarily acute,
and the midpoint $d$ of $AD$ does not lie on the Voronoi edge.
Instead, $d$ lies on the line of $vw$ so that $Avd$ has the same orientation
as its preimage under $\Gamma$, and $Adw$ has the opposite orientation.
Subtracting the integral over $Adw$ from the integral over $Avd$,
we get again the integral over $Avw$.
Adding the signed contributions of all triangles in the star,
we get the integral over Voronoi polygon.
In contrast to \eqref{eqn:detoverdet}, we formulate the claim by
integrating over all points of the \emph{image} of a triangle $\delta \in \Star{A}$.
We therefore replace the ratio of the determinants by its sign:
\begin{result}[Interior Cancellation Lemma]
  Let $A$ be an interior vertex in the Delaunay triangulation of a finite
  set $S \subseteq \Rspace^2$.
  Then
  $$
    \int_{x \in V_A} \dist{x}{A}^2 \diff x
       =  \sum_{\delta \in \Star{A}} \orientation{ \Gamma(\delta)}
          \int_{x \in \Gamma(\delta)} \dist{x}{A}^2 \diff x ,
  $$
  where $V_A$ is the Voronoi polygon of $A$,
  and the sum is over all triangles $\delta$ in the star of $A$
  inside the barycentric subdivision of the Delaunay triangulation.
\end{result}

To generalize the lemma to boundary vertices, we may clip the necessarily
infinite Voronoi polygons to within the convex hull of $S$.
With this modification, the lemma holds provided the angles opposite
to boundary edges are acute.
Indeed, the signed integrals over the triangles add up to the integral over
the cone from $A$ to the finite edges of the Voronoi polygon
plus the finite pieces of the infinite edges that end at the midpoints
of the boundary edges $GA$ and $AB$; see Figure \ref{fig:stars} on the right.
To summarize, we let $\N{S}{x}$ be the point in $S$ minimizing
the distance to $x$.
Let $\Dtri$ be the Delaunay triangulation of a finite set
$S \subseteq \Rspace^2$.
If all angles opposite to edges of the convex hull of $S$ are acute, then 
\begin{align}
  \Voronoif{\Dtri}  &=  \int_{x \in \conv{S}}
                       \dist{x}{\N{S}{x}}^2 \diff x .
  \label{eqn:SDT}
\end{align}
  
\paragraph{Difference of squared distances.}
To shed light on the general case, assume that $ABC$ is a triangle
with boundary edge $AB$ and obtuse angle at $C$.
The star of $A$ inside $\Sd{\Dtri}$ contains a linear sequence of triangles.
Removing the first triangle and the last, the signed areas add up
to the area of the cone of $A$ over the finite edges of its
Voronoi polygon.
Since the angle at $C$ is obtuse, the cone extends outside the convex hull
of $S$, but that extension is covered by the negatively oriented first
triangle in the linear sequence.
This triangle covers more, so we have a remaining negative contribution,
which is the integral over the region of points $x$ outside $\conv{S}$
for which the nearest vertex in $S$ is not the nearest vertex
on the boundary of $\conv{S}$; see Figure \ref{fig:stars} on the right.

To write this more succinctly, we recall that $\N{\Delta}{x}$
and $\NV{\Delta}{x}$ denote the nearest vertex and the nearest visible vertex
of $\Delta$ to $x \in \Rspace^2$.
Similarly, $\N{S}{x}$ is the nearest point in $S$ to $x$,
and we write $\NV{S}{x}$ for the nearest visible vertex of $\conv{S}$ to $x$.
For $x \in \conv{S}$, $\NV{S}{x}$ is not defined and we set
$\dist{x}{\NV{S}{x}}^2 = 0$ in this case.
\begin{result}[Voronoi Cell Decomposition Theorem]
  Let $\Dtri$ be the Delaunay triangulation of a finite set
  $S \subseteq \Rspace^2$.
  Then
  \begin{align}
    \Voronoif{\Dtri}  &=  \int_{x \in \Rspace^2}
       \left( \dist{x}{\N{S}{x}}^2 - \dist{x}{\NV{S}{x}}^2 \right) \diff x .
  \end{align}
\end{result}
\proof
 Write
 $g_{\Dtri} (x) = \dist{x}{\N{S}{x}}^2 - \dist{x}{\NV{S}{x}}^2$,
 and consider first the case in which all angles opposite to
 convex hull edges are acute.
 By \eqref{eqn:SDT},
 $\Voronoif{\Dtri}$ is the integral,
 over all points of $\conv{S}$, of $\dist{x}{\N{S}{x}}^2$.
 Since $\NV{S}{x}$ of a point in the convex hull is not defined,
 this is the same as the integral of
 $g_\Dtri (x)$, still only over the convex hull of $S$.
 For a point $x \not\in \conv{S}$,
 we have $\N{S}{x} = \NV{S}{x}$ by assumption on the angles.
 It follows that the contribution to the integral outside the convex hull
 vanishes, which implies the claimed equation in this special case.

 To extend the equation to the general case,
 we note that next to each convex hull edge but outside the convex hull,
 we get a region of points where the contributions to the nearest
 visible vertices do not cancel; see Figure \ref{fig:stars} on the right.
 Any two such regions are disjoint,
 and the net effect within each region is the described difference
 between the squared distances to the nearest point of $S$
 and the nearest visible vertex on the convex hull boundary.
 The claimed equation follows.
\eop

\paragraph{Optimality of the Delaunay triangulation.}
We prove that the Delaunay triangulation maximizes the Voronoi functional.
More than that, we show that the Delaunay triangulation maximizes
the functional locally, at every point of the plane.
To explain this, define
\begin{align}
  g_{\Delta} (x) &=  \dist{x}{\N{\Delta}{x}}^2 - \dist{x}{\NV{\Delta}{x}}^2 , \\
  g_{\Ktri} (x)  &=  \sum\nolimits_{\Delta \in \Ktri} g_{\Delta} (x) ,
\end{align}
where the sum is over all triangles $\Delta$ of a geometric triangulation
$\Ktri$ of $S$.
With this notation, we get
$\Voronoif{\Ktri} = \int_{x \in \Rspace^2} g_{\Ktri} (x) \diff x$
from \eqref{eqn:Voronoi-3}.
\begin{result}[Voronoi Optimality Theorem]
  Let $\Dtri$ be the Delaunay triangulation and $\Ktri$
  a geometric triangulation of a finite set $S \subseteq \Rspace^2$.
  Then $\Voronoif{\Ktri} \leq \Voronoif{\Dtri}$.
\end{result}
\proof
 We prove optimality by showing $g_{\Ktri} (x) \leq g_{\Dtri} (x)$
 for all $x \in \Rspace^2$.
 We have
 $g_{\Dtri} (x) = \dist{x}{\N{S}{x}}^2 - \dist{x}{\NV{S}{x}}^2$
 from the Voronoi Delaunay Theorem,
 so it suffices to show
 $g_{\Ktri} (x) \leq \dist{x}{\N{S}{x}}^2 - \dist{x}{\NV{S}{x}}^2$.
 Consider first a point $x \in \conv{S}$.
 For the triangle that contains $x$, we have
 $g_{\Delta} (x) = \dist{x}{\N{\Delta}{x}}^2$,
 which is positive.
 For every other triangle, we have $g_{\Delta} (x) \leq 0$.
 Let $A_0 = \N{S}{x}$, and consider the sequence of triangles
 $\Delta_0, \Delta_1, \ldots, \Delta_k$ constructed as follows.
 If $x$ lies inside one of the triangles in the star of $A_0$,
 then $\Delta_0$ is this triangle, and $k = 0$.
 Otherwise, there is a triangle $\Delta_0 \in \Star{A_0}$
 such that the line segment from $A_0$ to $x$ crosses one
 of its edges.
 Let $A_1$ be the nearest visible vertex of $\Delta_0$,
 and repeat the construction substituting $A_1$ for $A_0$ to get $\Delta_1$,
 and so on until $\Delta_k$ contains $x$ or we formed a cycle.
 In the former case, the sum of squared distances is
 \begin{align}
   \sum_{i=0}^k g_{\Delta_i} (x)  &\leq
       \sum_{i=0}^k \dist{x}{A_i}^2 - \sum_{i=0}^{k-1} \dist{x}{A_{i+1}}^2 .
 \end{align}
 The right hand side evaluates to $\dist{x}{A_0}^2 = \dist{x}{\N{S}{x}}^2$.
 The contribution of the triangles not in this sequence is non-positive,
 which implies
 $g_{\Ktri} (x) \leq \dist{x}{\N{S}{x}}^2 = g_{\Dtri} (x)$.
 If in the latter case the cycle goes around the triangle that contains $x$,
then it must contain a vertex, $A_\ell$, that is further from $x$
 than the vertices of the containing triangle, and in particular
 further than the closest vertex of that triangle, $y$.
 We stop the process at this vertex and get
 \begin{align}
   \sum_{i=0}^\ell g_{\Delta_i} (x)  &\leq
        \sum_{i=0}^{\ell-1} \dist{x}{A_i}^2
      - \sum_{i=0}^{\ell-1} \dist{x}{A_{i+1}}^2 .
 \end{align}
 This time, the right hand side evaluates to
 $\dist{x}{A_0}^2 - \dist{x}{A_\ell}^2 \leq \dist{x}{A_0}^2 - \dist{x}{y}^2$.
 The triangle that contains $x$ is the only one with a positive
 contribution, which is $\dist{x}{y}^2$, so that
 we again get $g_\Ktri (x) \leq \dist{x}{\N{S}{x}}^2$, as desired.

 \begin{figure}[hbt]
   \vspace{0.1in}
   \centering \resizebox{!}{1.4in}{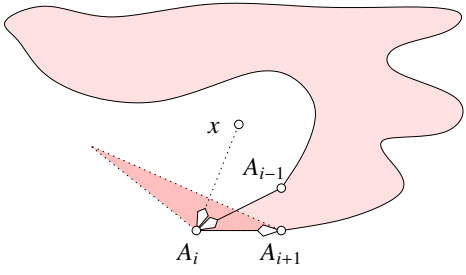}
   \caption{Assuming a cycle that does not go around $x$, there is a
     vertex $A_i$ such that both neighboring vertices lie on the same
     side of the line passing through $x$ and $A_i$.}
   \label{fig:Cycle}
 \end{figure}
 We still need to show that if the vertices $A_i$ form a cycle,
 then it goes around the point $x$.
 To derive a contradiction, we assume the cycle does not go around $x$,
 as in Figure \ref{fig:Cycle}.
 At each point $y$ of the cycle, we let $\varphi (y)$ be the
 counterclockwise rotation angle between the direction of the cycle at $y$
 and the vector $y - x$.
 This is a multivalued function, so we pick a branch by choosing
 its value at the starting point inside $[0, 2\pi)$
 and requiring that $\varphi$ be continuous along any edge,
 with right-continuous jumps at the vertices.
 The jump is positive at a counterclockwise turn and negative at a
 clockwise turn, never exceeding $\pi$ in absolute value.
 Assuming the cycle is oriented in a counterclockwise order,
 $\varphi$ grows by $2 \pi$ when we go around once.
 Indeed, because $x$ lies outside the cycle,
 the continuous changes total to $0$,
 and since the cycle is simple, the jumps total to $2 \pi$.
 Let $A_i$ be the vertex at which the jump changes $\varphi$ from
 a value less than $2 \pi$ to a value at least $2 \pi$.
 As illustrated in Figure \ref{fig:Cycle},
 the line passing through $A_{i-1}$ and $A_i$
 separates $x$ from $A_{i+1}$.
 But this implies that the edge connecting  $A_{i-1}$ to $A_i$
 crosses the edge of the triangle $\Delta_i$ opposite to $A_i$.
 This contradicts that both edges belong to $\Ktri$
 and concludes the proof that the cycle must go around $x$.

 To extend the argument to points $x$ outside the convex hull of $S$,
 we note that $g_{\Ktri} (x)$ and $g_{\Dtri} (x)$
 are both non-positive.
 Let $R$ be the region of points with $g_{\Dtri} (x) < 0$ and note
 that it is bounded.
 We can therefore add vertices and triangles on the outside,
 maintaining that we still have a Delaunay triangulation, such that
 $R$ is completely covered by the added triangles.
 Adding the same triangles to $\Ktri$,
 we get two new triangulations, $\Dtri'$ and $\Ktri'$,
 and we have $g_{\Ktri'} (x) \leq g_{\Dtri'} (x)$ for every $x \in R$
 using the above argument.
 Since we add the same triangles to $\Dtri$ and to $\Ktri$,
 we have
 $g{\Ktri'} (x) - g_{\Ktri} (x)
   = g_{\Dtri'} (x) - g_{\Dtri} (x)$,
 which implies $g_{\Ktri} (x) \leq g_{\Dtri} (x)$ also for the
 points $x \in R$.
 For the remaining points, $x \in \Rspace^2 \setminus \conv{S} \setminus R$,
 we have $g_{\Dtri} (x) = 0$ by definition of $R$,
 and $g_{\Ktri} (x) \leq 0$ because $x$ lies outside the underlying
 space of $\Ktri$.
 We thus get $g_{\Ktri} (x) \leq g_{\Dtri} (x)$
 for all points $x \in \Rspace^2$,
 as claimed.
\eop

\section{Non-optimality}
\label{sec4}

In this section, we show that the optimality of the Delaunay triangulation
among all geometric triangulation of a finite point set does not generalize
if we extend the family to all topological triangulations.
Furthermore, we show that even without this extension,
the optimality of the Delaunay triangulation does not extend to three dimensions.

\subsection{Topological Triangulations}
\label{sec41}

\begin{figure}[hbt]
  \centering \resizebox{!}{1.4in}{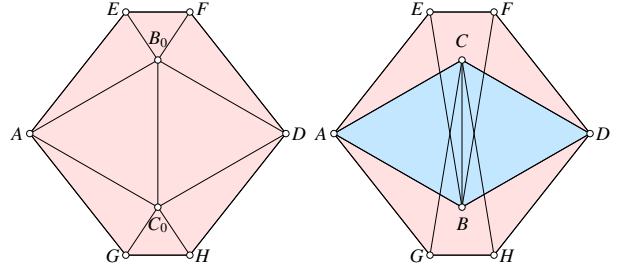}
  \caption{\emph{Left}: the Delaunay triangulation of the eight points.
    \emph{Right}: a topological triangulation of the same eight points
    obtained by exchanging the positions of points $B_0$ and $C_0$.}
  \label{fig:Folded}
\end{figure}

Consider the triangulation in Figure \ref{fig:Folded} on the left.
All triangles are acute, which implies that this is the Delaunay triangulation,
$\Dtri$, of the eight points.
The topological triangulation, $\Ktri$, of the same eight points is obtained
by exchanging the positions of points $B_0$ and $C_0$,
moving the edges and triangles along to preserve all incidences;
see Figure \ref{fig:Folded} on the right.
The pentagons $ACDFE$ and $ABDHG$ are covered only once by $\Ktri$,
while $ABDC$ is covered three times.
We prove $\Voronoif{\Ktri} > \Voronoif{\Dtri}$ by computing
$g_{\Ktri} (x) = \sum_{\Delta \in \Ktri} \orientation{\Delta} g_\Delta (x)$
for every point $x \in \Rspace^2$.
All triangles in $\Ktri$ are acute,
which implies that the only non-zero terms come from triangles that contain $x$.
First consider $x$ inside $ABC$, and note that the case $x$ inside $DBC$ is symmetric.
Here $x$ is covered by $ABC$, which has negative orientation,
and by two triangles with positive orientation,
namely one of $ABE$ and $BEF$ and one of $ACG$ and $CGH$.
In every case, $g_{\Ktri} (x)$ is the sum of the squared distances
to the nearest vertices of these triangles:
$g_{\Ktri} (x) = \alpha (x) + \beta (x) - \gamma (x)$, in which
\begin{align}
  \alpha (x)  &=  \min \left\{ \dist{x}{\N{ABE}{x}}^2 ,
                               \dist{x}{\N{BEF}{x}}^2 \right\} , \\
  \beta  (x)  &=  \min \left\{ \dist{x}{\N{ACG}{x}}^2 ,
                               \dist{x}{\N{CGH}{x}}^2 \right\} , \\
  \gamma (x)  &=  \dist{x}{\N{ABC}{x}}^2 .
\end{align}
We have $\alpha (x) \geq \dist{x}{\N{ABC}{x}}^2$ because $x$ is closer
to $C$ than to $E$ and $F$.
Similarly, $\beta (x) \geq \dist{x}{\N{ABC}{x}}^2$.
This implies
\begin{align}
  g_{\Ktri} (x)  &\geq  \dist{x}{\N{ABC}{x}}  =  g_{\Dtri} (x) ,
\end{align}
with strict inequality on a set of positive measure.
Next consider $x$ inside the region $ACDFE$,
and note that the case of $x$ inside $ABDHG$ is symmetric.
The region is covered only once, by triangles $ABE$, $EBF$, $DBF$,
all of which have positive orientation.
Since all three triangles are acute, $g_{\Ktri} (x)$ is the squared distance
of the closest vertex of the pentagon $ABDFE$,
while $g_{\Dtri} (x)$ is the squared distance to the closest vertex of $ACDFE$.
Point $B$ is further from all points within the pentagon than $C$,
which again implies $g_{\Ktri} (x) \geq g_{\Dtri} (x)$.
All triangles are acute, so $g_{\Ktri} (x) = g_{\Dtri} (x) = 0$
outside $AEFDHG$.
In summary, $g_{\Ktri} (x) \geq g_{\Dtri} (x)$ for all points $x \in \Rspace^2$,
and the inequality is strict on a set of positive measure.
Hence,
\begin{align}
  \Voronoif{\Ktri} &= \int_{\Rspace^2} g_{\Ktri} (x) \diff x > 
                     \int_{\Rspace^2} g_{\Dtri} (x) \diff x = \Voronoif{\Dtri},
\end{align}
as claimed.

\subsection{Beyond Two Dimensions}
\label{sec42}
We may generalize the Voronoi functional to three and higher dimensions,
defining it as the sum, over all simplices in the barycentric subdivision,
of the integral of the signed squared distance:
\begin{align}
  \Voronoif{\Ktri}  &=  \sum_{\delta \in \Sd{\Ktri}}
                       \orientation{\Gamma (\delta)}
                       \int_{x \in \Gamma (\delta)} \dist{x}{A}^2 \diff x , 
\end{align}
where $A$ is the unique vertex of $\delta$ that is also a vertex of $\Ktri$.
In $3$ dimensions, the circumcenter map moves two vertices of a
tetrahedron in the barycentric subdivision, possibly inverting
a tetrahedron twice and thus returning it to positive orientation.
This explains why there may be triangulations for which the
Voronoi functional exceeds that of the Delaunay triangulation.

\begin{figure}[hbt]
  \centering \resizebox{!}{2.2in}{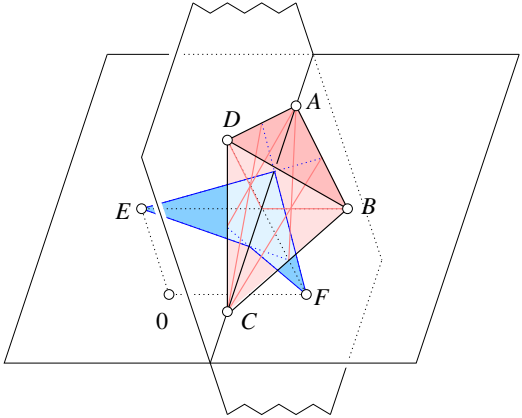}
  \caption{A tetrahedron with partial barycentric subdivision and partial
    image under the circumcenter map.
    Two of the triangles have obtuse face angles, which we highlight
    by showing the containing planes and the centers of the circumcircles.}
  \label{fig:BadTet}
\end{figure}

\paragraph{Double fold-over.}
We exhibit a tetrahedron for which there are points not contained in the
tetrahedron that have a positive contribution to the integral
of the squared distance to the nearest vertex.
In $\Rspace^2$, a triangle with this property does not exist.
Figure \ref{fig:BadTet} shows the tetrahedron, $ABCD$, together with
the center of the circumsphere, $0 = \circumcenter{ABCD}$.
The faces $ABD$, $CBD$ are isosceles triangles with three acute angles each.
In contrast, $BAC$, $DAC$ are isosceles triangles with obtuse angles at $B$ and at $D$,
and we show the centers of their circumcircles:
$E = \circumcenter{BAC}$ and $F = \circumcenter{DAC}$.
The circumcenter map leaves the four vertices and the midpoints of
the edges fixed,
and it moves the barycenters of the four triangles and
of the tetrahedron to new locations.
As a consequence, some of the images of the $24$ triangles in the
barycentric subdivision of the boundary of $ABCD$ change their
orientation within their respective planes.
Because of the obtuse angles at $B$ and at $D$,
these are the four triangles that share the midpoint of $AC$.
The other $20$ triangles preserve their orientation.
The image of the barycenter of the tetrahedron is $0$,
the planes of $BAC$ and $DAC$ separate $\barycenter{ABCD}$ from $0$,
while the planes of $ABD$ and $CBD$ do not separate the two points.
It follows that $16$ tetrahedra in the barycentric subdivision
preserve their orientation, while $8$ tetrahedra reverse their orientation.
Just like in $\Rspace^2$, where we have regions outside an obtuse triangle
whose points contribute negatively to the Voronoi functional
(see Figure \ref{fig:fold}),
we have regions outside $ABCD$ whose points contribute positively
to the Voronoi functional.
Specifically, there is such a point $x$ near $E$ inside the cone of $0$
over the blue triangle in Figure \ref{fig:BadTet}.
It satisfies $\dist{x}{B} < \dist{x}{C} < \dist{x}{A}$,
and its contribution to the Voronoi functional is
$\dist{x}{C}^2 - \dist{x}{B}^2$, which is positive.
The argument used to prove the optimality of the Delaunay triangulation
in $\Rspace^2$ does therefore not apply.

\begin{figure}[hbt]
  \vspace{-0.2in}
  \hspace{-0.10in} \includegraphics[width=0.45\textwidth]{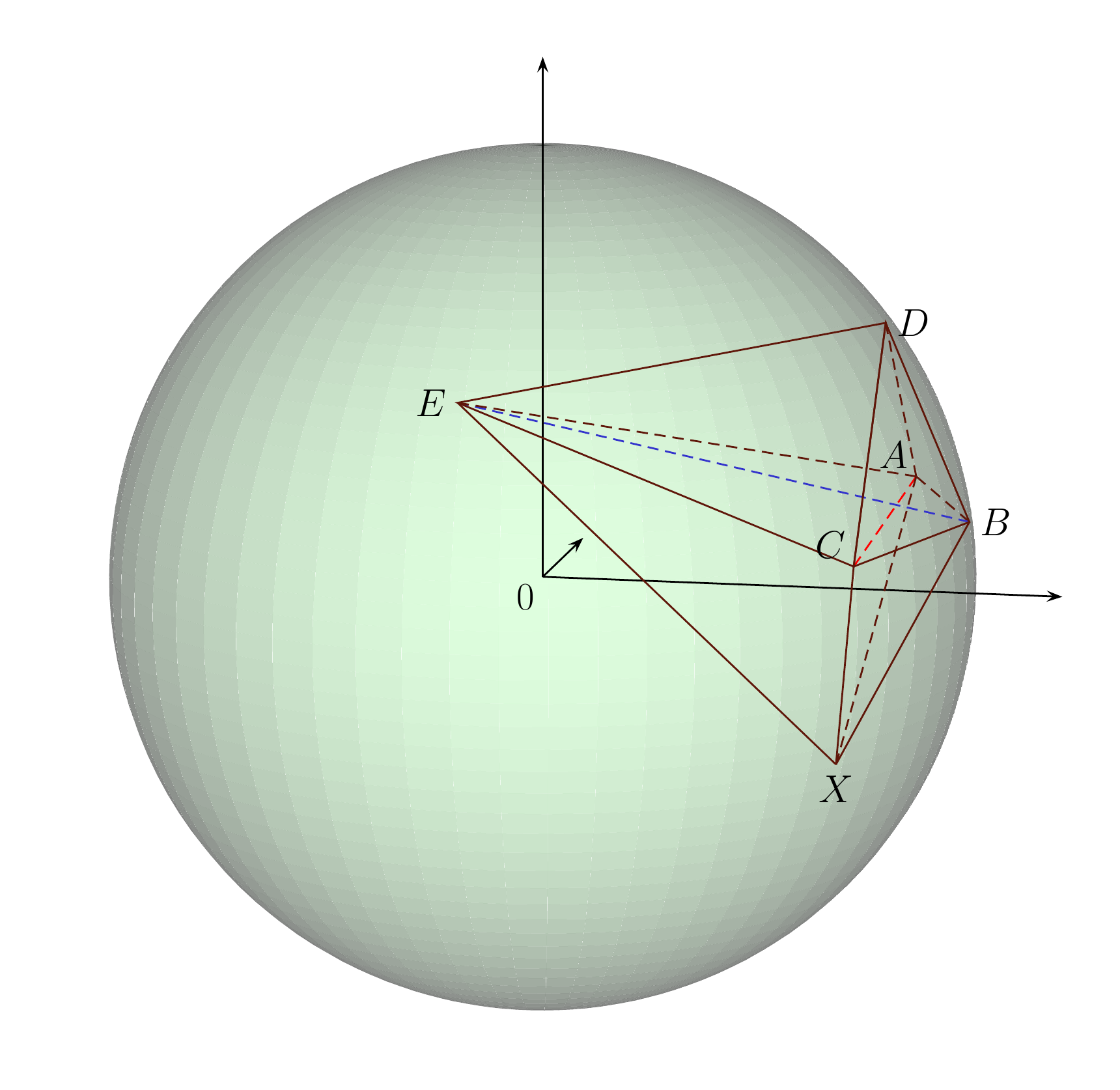}
  \vspace{-0.2in}
  \caption{The six points are the vertices of a non-regular octahedron,
    which we decompose into four tetrahedra in two ways.}
  \label{fig:counter}
\end{figure}

\paragraph{Counterexample.}
Inspired by the construction in Figure \ref{fig:BadTet},
we use a numerical method for estimating the Voronoi functional
to find points in $\Rspace^3$ for which the Delaunay triangulation
does not maximize the functional.
Figure \ref{fig:counter} illustrates one such example.
There are six point with coordinates
\begin{align*}
 \begin{array}{rrrrrrr}
  A &=& (\!\!\!\!\!\!  &   7.99,  & \hspace{-0.1in}  5.80,
                       & \hspace{-0.1in}  1.65 & \hspace{-0.12in}),  \\
  B &=& (\!\!\!\!\!\!  &   9.86,  & \hspace{-0.1in}  0.00,
                       & \hspace{-0.1in}  1.65 & \hspace{-0.12in}),  \\
  C &=& (\!\!\!\!\!\!  &   7.80,  & \hspace{-0.1in} -5.80,
                       & \hspace{-0.1in}  1.65 & \hspace{-0.12in}),  \\
  D &=& (\!\!\!\!\!\!  &   7.89,  & \hspace{-0.1in}  0.00, 
                       & \hspace{-0.1in}  6.14 & \hspace{-0.12in}),  \\
  E &=& (\!\!\!\!\!\!  &  -2.00,  & \hspace{-0.1in} -0.01,
                       & \hspace{-0.1in}  4.02 & \hspace{-0.12in}),  \\
  X &=& (\!\!\!\!\!\!  &   6.89,  & \hspace{-0.1in}  0.00,
                       & \hspace{-0.1in} -4.14 & \hspace{-0.12in}),
 \end{array}
\end{align*}
whose convex hull is an octahedron.
The Delaunay triangulation decomposes the octahedron into the four
tetrahedra sharing the edge $BE$.
We compare it with the triangulation, $\Ktri$,
whose four tetrahedra share the edge $AC$.
The values of the Voronoi functional are
\begin{align}
  \Voronoif{\Ktri}  &=  3432.96 \pm 0.01 , \\
  \Voronoif{\Dtri}  &=  3413.75 \pm 0.01 ,
\end{align}
which shows that $\Dtri$ does not maximize the functional.

\section{Discussion}
\label{sec5}

As a by-product of the Voronoi Optimality Theorem, we have the optimality of
the Delaunay triangulation for the functional that takes the sum over all
triangles of the squared circumradius times the area.
More generally, we define
\begin{align}
  \Radiusf{\alpha}(\Ktri)  &=  \sum_{\Delta \in \Ktri} R_\Delta^\alpha \area{\Delta} ,
\end{align}
where $R_\Delta$ is the circumradius of $\Delta$.
Indeed, $\Radiusf{2}(\Ktri) = 3 \Rajanf{\Ktri} - 3 \Voronoif{\Ktri}$
so we conclude that $\Radiusf{2}$ attains its minimum at the Delaunay triangulation.
Similarly, we can prove that $\Radiusf{1}$ attains its
minimum at the Delaunay triangulation.
Writing $p(ABC) = \dist{A}{B} \cdot \dist{B}{C} \cdot \dist{C}{A}$,
we have $4 \area{\Delta} = p(\Delta)/R_\Delta$ and therefore
$\Radiusf{1}(\Ktri) = \frac{1}{4} \sum_{\Delta \in \Ktri} p(\Delta)$.
When we flip an edge to turn a geometric triangulation into the Delaunay triangulation,
the sum of $p(\Delta)$ cannot increase,
which implies the Delaunay triangulation is optimal.
It would be interesting to prove optimality for all $\alpha \geq 1$.
\ignore{
Finally, we consider
\begin{align}
  \Delaunayf{\Ktri}  &=  \sum_{\Delta \in \Ktri}
    \dist{\barycenter{\Delta}}{\circumcenter{\Delta}}^2 \area{\Delta} ,
\end{align}
where $\barycenter{\Delta}$ is the barycenter and $\circumcenter{\Delta}$
is the circumcenter of the triangle.
Again we have computational evidence that it attains its minimum at
the Delaunay triangulation.
It would be interesting to prove that this is the case.
}

It might be interesting to generalize the Voronoi functional to points
with weights; see e.g.\ \cite{Ede01}.
The inclusion-exclusion argument presented in Section \ref{sec3}
extends to the weighted case, but is it true that the generalized
Voronoi functional is maximized by the Delaunay triangulation
for weighted points in $\Rspace^2$?

\newpage

\clearpage
\appendix
\section{Optimality by Flipping}
\label{appA}

\begin{figure}[hbt]
   \centering \resizebox{!}{2.6in}{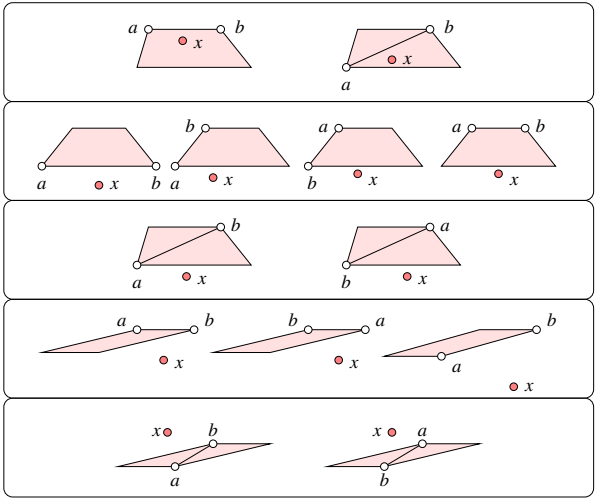}
   \caption{The cases considered by the flip algorithm.
     From \emph{top} to \emph{bottom}: Case I, Case II.1, Case II.2,
     Case III.1, Case III.2.
     We draw the diagonal whenever we know that it belongs to the
     Delaunay triangulation.}
   \label{fig:Cases}
 \end{figure}

As an alternative to the proof in Section \ref{sec3},
we can show the optimality of Delaunay triangulations
for the Voronoi functional via flips.
More precisely, it is sufficient to show that the Voronoi Optimality Theorem
holds for any convex quadrangle.
Since the Delaunay triangulation may be derived from any triangulation
through a sequence of flips,
this will imply that the theorem is true for all finite sets
$S \subseteq \Rspace^2$.
\begin{result}[Voronoi Optimality Theorem for Flips]
  Let $\Dtri$ and $\Ktri$ be the Delaunay and the non-Delaunay
  triangulations of four points in convex position in $\Rspace^2$.
  Then $g_{\Ktri}(x) \leq g_{\Dtri}(x)$ for all $x \in \Rspace^2$.
\end{result}
\proof
 Let $a, b, c, d \in \Rspace^2$ be the four points in convex position,
 $x \in \Rspace^2$, and assume
 $\dist{x}{a} \leq \dist{x}{b} \leq \dist{x}{c} \leq \dist{x}{d}$.
 Then
 \begin{enumerate}\denselist
   \item[{[A]}] the edge $ab$ belongs to the Delaunay triangulation,
     because there is a circle (centered at $x$) such that $a, b$
     are on or inside, and $c, d$ are on or outside the circle,
   \item[{[B]}] the points $c$ and $d$ do not lie inside the triangle $abx$,
     else they would be closer to $x$ than $a$ or $b$.
 \end{enumerate}
 [A] implies that if $ab$ is a diagonal of the quadrangle,
 then $\Dtri$ consists of the triangles $abc, abd$.
 The remainder of the proof is a case analysis, which we illustrate
 in Figure \ref{fig:Cases}.
 \begin{description}\denselist
   \item[{\sc Case I}] $x$ is inside the quadrangle;
     see Figure \ref{fig:Cases}, $1$st row.
     \begin{description}\denselist
       \item[{\sc Case I.1}] $ab$ is a side of the quadrangle.
         Then $g_\Dtri (x) = g_\Ktri (x)$ because the values for the
         two triangles are either $\dist{x}{a}^2$ and $0$,
         or $\dist{x}{b}^2$ and $\dist{x}{a}^2-\dist{x}{b}^2$.
       \item[{\sc Case I.2.}] $ab$ is a diagonal.
         Using [A], we get $g_{\Dtri}(x) = \dist{x}{a}^2$,
         which is greater than or equal to
         $g_\Ktri (x) = \dist{x}{a}^2 + \dist{x}{b}^2-\dist{x}{c}^2$.
       \end{description}

    \item[{\sc Case II}]  $x$ lies outside,
      and only two vertices of the quadrangle are visible from $x$.
      \begin{description}\denselist
        \item[{\sc Case II.1}]  $ab$ is a side;
          see Figure \ref{fig:Cases}, $2$nd row.
          We have four subcases and $g_\Dtri (x) = g_\Ktri (x)$ in each.
          \begin{description}\denselist
            \item[{\sc Case II.1.a}]  Both $a$ and $b$ are visible from $x$.
              Then $g_\Delta (x)=0$ for all four triangles.
            \item[{\sc Case II.1.b}]  $a$ is visible and $b$ is invisible.
              Then $g_\Delta (x)=0$ for all four triangles.
            \item[{\sc Case II.1.c}]  $b$ is visible and $a$ is invisible.
              In both triangulations, the values for the two triangles
              are $0$ and $\dist{x}{a}^2 - \dist{x}{b}^2$.
            \item[{\sc Case II.1.d}]  Both $a$ and $b$ are invisible.
              For one triangulation the values for the two triangles are
              $\dist{x}{a}^2 - \dist{x}{b}^2$ and
              $\dist{x}{b}^2 - \dist{x}{c}^2$,
              while for the other triangulation the values are
              $\dist{x}{a}^2-\dist{x}{c}^2$ and $0$.
          \end{description}

        \item[{\sc Case II.2}]  $ab$ is a diagonal;
          see Figure \ref{fig:Cases}, $3$rd row.
          We have two subcases and
          $g_\Dtri (x) = g_\Ktri (x) - \dist{x}{b}^2 + \dist{x}{c}^2$ in each.
          By [A] we know which triangulation is Delaunay.
          \begin{description}\denselist
            \item[{\sc Case II.2.a}]  $a$ is visible and $b$ is invisible.
              Then $g_{\Dtri}(x) = 0$ and
              $g_\Ktri (x) = \dist{x}{b}^2 - \dist{x}{c}^2$.
            \item[{\sc Case II.2.b}]  $b$ is visible and $a$ is invisible.
              Then $g_{\Dtri}(x) = \dist{x}{a}^2-\dist{x}{b}^2$ and
              $g_\Ktri (x) = \dist{x}{a}^2-\dist{x}{c}^2$.
          \end{description}
      \end{description}

    \item[{\sc Case III}]  $x$ is outside the quadrangle and three vertices
      of the quadrangle are visible from $x$.
      \begin{description}\denselist
        \item[{\sc Case III.1}]  $ab$ is a side;
          see Figure \ref{fig:Cases}, $4$th row.
          We have three subcases and $g_\Dtri (x) = g_\Ktri (x)$ in each.
          \begin{description}\denselist
            \item[{\sc Case III.1.a}]  $a$ is invisible.
              By [B] there is only one choice for $b$
              so that the line through $x$ and $a$ separates $b$ from $c, d$.
              Both triangulations contain one triangle with
              value $0$ and the other with $\dist{x}{a}^2-\dist{x}{b}^2$.
            \item[{\sc Case III.1.b}]  $b$ is invisible.
              As in Case III.1.a, there is only one choice for $a$,
              and $g_\Delta(x) = 0$ for all four triangles.
            \item[{\sc case III.1.c}]  $a$ and $b$ are both visible.
              Then $g_\Delta(x)=0$ for all four triangles.
          \end{description}
        \item[{\sc Case III.2}]  $ab$ is a diagonal;
          see Figure \ref{fig:Cases}, $5$th row.
          We have two subcases and 
          $g_\Dtri (x) = g_\Ktri (x) - \dist{x}{b}^2 + \dist{x}{c}^2$ in each.
          By [A] we know which triangulation is Delaunay.
          \begin{description}\denselist
            \item[{\sc Case III.2.a}]  $a$ is invisible.
              $g_{\Dtri}(x) = \dist{x}{a}^2 - \dist{x}{b}^2$ and
              $g_\Ktri (x) = \dist{x}{a}^2 - \dist{x}{c}^2$.
            \item[{\sc Case III.2.b}]  $b$ is invisible.
              Then $g_{\Dtri}(x) = 0$ and
              $g_\Ktri (x) = \dist{x}{b}^2-\dist{x}{c}^2$.
            \item[{\sc Case III.2.c}]  That $a$ and $b$ are both visible
              is impossible because it would contradict [B].
          \end{description}
      \end{description}
 \end{description}
 The case analysis is exhaustive, which implies
 $g_{\Dtri}(x) \geq g_{\Ktri}(x)$, as claimed.
\eop

We remark that a flip can have only two possible outcomes:
it leaves the value of $g_\Ktri (x)$ the same,
or it increases this value by $- \dist{x}{b}^2 + \dist{x}{c}^2$.
The latter case happens iff the two vertices
closest to $x$ form a diagonal.

\section{Calculations}
\label{appB}

In this section, we derive formula \eqref{eqn:Voronoi-2},
both for the acute and obtuse cases.
We will be working in $\Rspace^3$ and denote the three Cartesian coordinates
by $x, y, z$.
The Voronoi functional is invariant under planar isometries,
so we may assume that the circumcenter of the triangle is at the origin.
The points $A', B', C'$ thus have a common $z$-coordinate, which is $R^2$,
and the $z$-coordinate of $0'$ is $-R^2$.
We denote the angles of $ABC$ by $\alpha, \beta, \gamma$
and the area of $ABC$ by $\sigma$.
We begin with the formula from Section \ref{sec2}:
\begin{align}
  \Voronoif{ABC} &= \volume{ABC} - \Rajanf{ABC} \nonumber \\
                   &- \volume{Q_A} - \volume{Q_B} - \volume{Q_C}
  \label{eqn:vor_acute}
\end{align}
We divide the quadrangle $Q_A$ into two triangles, $AB_10$ and $AC_10$,
and similarly for $Q_B$ and $Q_C$.
The triangles adjacent to side $AB$ are congruent,
and we denote both of them by $T_{AB}$.
Similarly, two of the other triangles are denoted by $T_{BC}$, and two by $T_{AC}$.
Now we rewrite \eqref{eqn:vor_acute}:
\begin{align}
  \Voronoif{ABC} &= \volume{ABC} - \Rajanf{ABC}  \nonumber \\
                   &- 2\volume{T_{AB}} - 2\volume{T_{BC}} - 2\volume{T_{AC}}
  \label{eqn:vor_acute2}
\end{align}
Since $ABCA'B'C'$ is a prism, $\volume{ABC} = \sigma R^2$.
To compute $\volume{T_{AB}}$, $\volume{T_{BC}}$, $\volume{T_{AC}}$,
we use the following result.
\begin{result}[Angle Lemma]
  Let $KLM$ be a triangle with vertices $K=(0,0,-R^2)$, $L=(R,0,R^2)$,
  and the plane $KLM$ perpendicular to the $xz$-plane.
  Then we have $\angle KML = \pi/2$
  $\volume{KLM} = \tfrac{1}{24} R^4\sin{4\phi}$, in which $\phi = \angle MKL$.
\end{result}
\proof
 We draw the altitude $MH$ to divide $KLM$ into $KHM$ and $LHM$.
 Then we can find the area under this triangle for each plane perpendicular
 to the $x$-axis and integrate it:
 \begin{align*}
   \volume{KHM} &= \int_{0}^{R \cos^2{\phi}} x \tan{\phi} (2Rx-R^2) \diff x   \\
                &= \left[ \tan{\phi} \left( \tfrac{2}{3} R x^3
                                   - \tfrac{1}{2} R^2 x^2 \right)
                                     \right]_0^{R\cos^2{\phi}}                \\
                &= R^4 \tan{\phi} \left( \tfrac{2}{3} \cos^6{\phi}
                                   - \tfrac{1}{2} \cos^4{\phi} \right).
 \end{align*}
 Similarly,
 $\volume{LHM} = - R^4 \cot{\phi} \left( \tfrac{2}{3} \sin^6{\phi}
                            - \tfrac{1}{2} \sin^4{\phi} \right)$.
 Writing $\kappa = \volume{KLM} = \volume{LHM} + \volume{KHM}$, we get
 \begin{align*}
   \kappa &= R^4 \tan{\phi} \left( \tfrac{2}{3} \cos^6{\phi}
           - \tfrac{1}{2} \cos^4{\phi} \right)                      
           - R^4 \cot{\phi} \left( \tfrac{2}{3} \sin^6{\phi}
           - \tfrac{1}{2} \sin^4{\phi} \right)                      \\
          &= R^4 \sin{\phi} \cos{\phi} \left( \tfrac{2}{3} \cos^4{\phi}
           - \tfrac{1}{2} \cos^2{\phi}                       
           - \tfrac{2}{3} \sin^4{\phi}
           + \tfrac{1}{2} \sin^2{\phi} \right)                      \\
          &= R^4 \sin{\phi} \cos{\phi}
             \left( \tfrac{2}{3} (\cos^4{\phi} - \sin^4{\phi}) 
           - \tfrac{1}{2} (\cos^2{\phi}
           - \sin^2{\phi}) \right)                                  \\
          &= \tfrac{1}{6} R^4 \sin{\phi} \cos{\phi}
             \left( \cos^2{\phi} - \sin^2{\phi}\right)              \\
          &= \tfrac{1}{12} R^4 \sin{2\phi} \cos{2\phi}             
 \end{align*}
 which evaluates to $\kappa = \tfrac{1}{24} R^4 \sin{4\phi}$, as claimed.
\eop

Writing $\upsilon = 2\volume{T_{AB}} + 2\volume{T_{BC}} + 2\volume{T_{AC}}$,
we apply the lemma to all six triangles to obtain
\begin{align*}
  \upsilon  &= \tfrac{1}{12} R^4 (\sin{4\alpha} + \sin{4\beta} + \sin{4\gamma}) \\
     &= \tfrac{1}{12} R^4 (2\sin{(2\alpha+2\beta)} \cos{(2\alpha-2\beta)}
         + \sin{4\gamma})                                                \\
     &= \tfrac{1}{12} R^4 (-2\sin{2\gamma} \cos{(2\alpha-2\beta)}
         + 2\sin{2\gamma} \cos{2\gamma})                                 \\
     &= \tfrac{1}{3} R^4 \sin{2\gamma} \sin{(\alpha-\beta-\gamma)}
        \sin{(\alpha+\gamma - \beta)}                                    \\
     &= -\tfrac 8 3 R^4 \sin{\alpha} \sin{\beta} \sin{\gamma}
        \cos{\alpha} \cos{\beta} \cos{\gamma}                            \\
     &= -\tfrac{1}{3} Rabc \cos{\alpha} \cos{\beta} \cos{\gamma}           \\
     &= -\tfrac{4}{3} R^2\sigma \cos{\alpha} \cos{\beta} \cos{\gamma}.
\end{align*}
Finally, we combine it with \eqref{eqn:vor_acute2} and \eqref{rajan}:
\begin{align*}
  Vf(ABC) &= \sigma R^2 - \tfrac{1}{12}
             \sigma R^2 (4\sin^2{\alpha} + 4\sin^2{\beta}                \\
          &+ 4\sin^2{\gamma})+ \tfrac{4}{3} \sigma R^2 \cos{\alpha} \cos{\beta} \cos{\gamma}.
\end{align*}
We use the identities
\begin{align*}
  \cos{\alpha} \cos{\beta} \cos{\gamma} &= \frac{p^2-(2R+r)^2}{4R^2} , \\
  \sin^2{\alpha} + \sin^2{\beta} + \sin^2{\gamma}
                 &= \frac{p^2-r^2-4Rr}{2R^2},
\end{align*}
where $p$ is the half-perimeter of $ABC$, and $r$ is its inradius
(for details see \cite[12.43a,b]{Pra06}).
\begin{align*}
  Vf(ABC) &= \sigma \left( R^2 + \tfrac{1}{3} (p^2-(2R+r)^2)
           - \tfrac{1}{6} (p^2-r^2-4Rr)\right)                             \\
          &= \sigma \left( \tfrac{1}{6} (p^2-r^2-4Rr)
           - \tfrac{1}{3} R^2\right)                                       \\
          &= \tfrac{1}{12} \sigma (a^2+b^2+c^2) - \tfrac{1}{3} \sigma R^2.
\end{align*}
The sequence of calculations for the obtuse triangle is almost the same.
Assume that $\angle ABC > \pi/2$.
Then the formula \eqref{eqn:vor_acute2} is transformed to:
\begin{align}
  \Voronoif{ABC} &= \volume{ABC} - \Rajanf{ABC}       \nonumber       \\
                   &- 2\volume{T_{AB}} - 2\volume{T_{BC}} + 2\volume{T_{AC}} .
  \label{eqn:vor_obtuse}
\end{align}
In this case, $\volume{T_{AC}} = \tfrac{1}{24} R^4 \sin{4(\pi-\beta)}
                               = -\tfrac{1}{24} R^4 \sin{4\beta}$.
Therefore, all further calculations are completely the same as in the acute case.

\newpage
\section{Notation}
\label{appC}

\begin{table}[hbt]
  \centering
  \begin{tabular}{ll}
    $S \subseteq \Rspace^2$
      &  finite set of points                                     \\
    $\Voronoif{\Delta}$
      &  Voronoi functional of triangle                           \\
    $\varpi \colon \Rspace^2 \to \Rspace$
      &  unit paraboloid                                          \\
    $f_A \colon \Rspace^2 \to \Rspace$
      &  tangent plane                                            \\
    $A, B, C$
      &  vertices                                                 \\
    $A_1, B_1, C_1, 0$
      &  barycenters                                              \\
    $A', A_1', 0'$
      &  lifted points                                            \\
                                                                  \\
    $\ksx \in \Ktri$
      &  simplex in triangulation                                 \\
    $\delta \in \Sd{\Ktri}$
      &  triangle in barycentric subdivision                      \\
    $\barycenter{\ksx}, \circumcenter{\ksx}, \circumradius{\ksx}$
      &  barycenter, circumcenter, circumradius                   \\
    $\Gamma \colon \us{\Sd{\Ktri}} \to \Rspace^2$
      &  circumcenter map                                         \\
    $\Height \colon \us{\Sd{\Ktri}} \to \Rspace$
      &  height map                                               \\
                                                                  \\
    $\N{\Delta}{x}, \NV{\Delta}{x}$
      &  nearest, nearest visible vertex                          \\
    $\N{S}{x}, \NV{S}{x}$
      &  nearest, nearest visible point                           \\
    $\Dtri$
      &  Delaunay triangulation                                   \\
    $\Delta_0, \ldots, \Delta_k$
      &  sequence of triangles                                    \\
    $g_{\Delta} (x), g_{\Ktri} (x)$
      &  local Voronoi functional
  \end{tabular}
  \caption{Notation for geometric concepts, sets, functions,
           vectors, variables used in the paper.}
  \label{tbl:Notation}
\end{table}

\newpage
\section{Definitions and Claims}
\label{appD}

\begin{itemize}\denselist
  \item  Section \ref{sec1}:  Introduction.
  \item  Section \ref{sec2}:  The Voronoi Functional.
    \begin{itemize}\denselist
      \item  $\Voronoif{\Delta}$.
      \item  $\Rajanf{\Delta}$.
      \item  $\Voronoif{\Ktri}$.
    \end{itemize}
  \item  Section \ref{sec3}:  Geometric Triangulations.
    \begin{itemize}\denselist
      \item  \emph{Circumcenter map}.
      \item  \emph{Height map}.
      \item  {\sc Interior Cancellation Lemma}.
      \item  {\sc Voronoi Cell Decomposition Theorem}.
      \item  {\sc Voronoi Optimality Theorem}.
    \end{itemize}
  \item  Section \ref{sec4}:  Non-optimality.
    \begin{itemize}\denselist
      \item  Section \ref{sec41}:  Topological Triangulations.
      \item  Section \ref{sec42}:  Beyond Two Dimensions.
    \end{itemize}
  \item  Section \ref{sec5}:  Discussion.
    \begin{itemize}\denselist
      \item  $\Radiusf{\alpha}(\Delta)$.
    \end{itemize}
  \item  Section \ref{appA}:  Voronoi Optimality by Flipping.
    \begin{itemize}\denselist
      \item  {\sc Voronoi Optimality Theorem for Flips}.
    \end{itemize}
  \item  Section \ref{appB}:  Calculations.
    \begin{itemize}\denselist
      \item  {\sc Angle Lemma}.
    \end{itemize}
\end{itemize}


\end{document}